\documentclass[11pt]{amsart}

\usepackage{amsmath,latexsym,amssymb,amsthm,array,amsfonts}

\usepackage{mathrsfs,dsfont,bbm}

\theoremstyle{plain}

\numberwithin{equation}{section}

\theoremstyle{plain}

\newtheorem{theorem}{Theorem}[section]
\newtheorem{definition}{Definition}[section]
\newtheorem{lemma}{Lemma}[section]
\newtheorem{corollary}{Corollary}[section]
\newtheorem{proposition}{Proposition}[section]

\newtheorem*{acknowledgement}{Acknowledgement}

\begin{document}
\title[Integration with respect to fractional local time]
{Integration with respect to fractional local time with Hurst index
$1/2<H<1$${}^{*}$}

\footnote[0]{${}^{*}$The Project-sponsored by NSFC (10571025) and
the Key Project of Chinese Ministry of Education (No.106076).}

\author[L. Yan, J. Liu and X. Yang]
{Litan Yan${}^{1,\dag}$, Junfeng Liu${}^{2}$ and Xiangfeng
Yang${}^{1}$}

\footnote[0]{${}^{\dag}$Corresponding Author (litanyan@hotmail.com)}

\date{}



\maketitle

\date{}

\begin{center}
{\it ${}^1$Department of Mathematics, College of Science, Donghua
University\\
2999 North Renmin Rd., Songjiang, Shanghai 201620, P. R. China}\\
{\it ${}^2$Department of Mathematics, East China University of
Science and Technology\\
130 Mei Long Rd., Xuhui, Shanghai 200237, P.R. China}
\end{center}

\maketitle

\begin{abstract}
Let ${\mathscr L}^H(x,t)=2H\int_0^t\delta(B^H_s-x)s^{2H-1}ds$ be the
weighted local time of fractional Brownian motion $B^H$ with Hurst
index $1/2<H<1$. In this paper, we use Young integration to study
the integral of determinate functions $\int_{\mathbb R}f(x){\mathscr
L}^H(dx,t)$. As an application, we investigate the {\it weighted
quadratic covariation} $[f(B^H),B^H]^{(W)}$ defined by
$$
[f(B^H),B^H]^{(W)}_t:=\lim_{n\to \infty}2H\sum_{k=0}^{n-1}
k^{2H-1}\{f(B^H_{t_{k+1}})-f(B^H_{t_{k}})\}
(B^H_{t_{k+1}}-B^H_{t_{k}}),
$$
where the limit is uniform in probability and $t_k=kt/n$. We show
that it exists and
$$
[f(B^H),B^H]^{(W)}_t=-\int_{\mathbb R}f(x){\mathscr L}^H(dx,t),
$$
provided $f$ is of bounded $p$-variation with $1\leq
p<\frac{2H}{1-H}$. Moreover, we extend this result to the
time-dependent case. These allow us to write the fractional It\^{o}
formula for new classes of functions.
\end{abstract}

\section{Introduction}
The fractional It\^o's formula has played a central role in
stochastic analysis for fractional Brownian motion and almost all
aspects of its applications. But, the restriction of It\^o's formula
to functions with twice differentiability often encounter
difficulties in applications. Extensions to less smooth functions
are useful in studying many problems such as partial differential
equations with some singularities and mathematics of finance. For
the classical It\^o's formula, many authors have written extended
versions of the It\^o formula in order to relax this regularity
condition. Still there are always some new problems requiring the
use of "It\^o formula" under lighter conditions. One approach of
extending It\^{o}'s formula is by using local time-space calculus to
absolutely continuous function $F$ with locally bounded measurable
derivative $F'$ (see Bouleau--Yor~\cite{Boul},
Eisenbaum~\cite{Eisen1,Eisen2}, Feng-Zhao~\cite{Feng}, F\"ollmer
{\it et al}~\cite{Follmer}, Peskir~\cite{Peskir1},
Russo--Vallois~\cite{Russo2}, Yan--Yang~\cite{Yan2}, and the
references therein). Moreover, the backward integral and quadratic
covariation are fundamental tools in these discussions. However, all
these extensions are not effective on fractional Brownian motion
$B^H$ since the quadratic covariation $[f(B^H),B^H]$ of $f(B^H)$ and
$B^H$ satisfies
$$
\left[f(B^{H}),B^{H}\right]_t=
\begin{cases}
0, &{\text {if $\frac12<H<1$}}\\
+\infty, &{\text {if $0<H<\frac12$}}
\end{cases}
$$
for all $t\geq 0$, in general, where $x\mapsto f(x)$ is a
determinate function. As well-known, the fractional Brownian motion
(fBm) with Hurst index $H\in (0,1)$ is a mean zero Gaussian process
$B^H=\{B_t^H, 0\leq t\leq T\}$ with the covariance function
$$
E\left[B_t^HB_s^H\right]=\frac{1}{2}\left[t^{2H}+s^{2H}-|t-s|^{2H}
\right]
$$
for all $t,s\geqslant 0.$ For $H=1/2$, $B^H$ coincides with the
standard Brownian motion $B$. $B^H$ is neither a semimartingale nor
a Markov process unless $H=1/2$, so many of the powerful techniques
from stochastic analysis are not available when dealing with $B^H$.
However, as a Gaussian process, one have constructed the stochastic
calculus of variations with respect to $B^{H}$ (see Biagini {\it et
al}~\cite{BHOZ} and references therein).

In order to motivate our subject, let us first recall some known
results concerning the quadratic variation and It\^o's formula. Let
$F$ be an absolutely continuous function with locally square
integrable derivative $f$, that is,
$$
F(x)=F(0)+\int_0^xf(y)dy
$$
with $f$ being locally square integrable. F\"ollmer {\it et
al}~\cite{Follmer} introduced the following It\^o's formula:
\begin{equation}\label{sec1-eq1.1}
F(B_t)=F(0)+\int_0^tf(B_s)dB_s+\frac12\left[f(B),B\right]_t.
\end{equation}
If $f\in C^1({\mathbb R})$,~\eqref{sec1-eq1.1} is the classical
It\^o formula. Moreover, the result has been extended to some
semimartingales and smooth nondegenerate martingales (see
Russo--Vallois~\cite{Russo2} and Moret--Nualart~\cite{Moret}).
However, it is obvious that the formula~\eqref{sec1-eq1.1} is not
true for fBm $B^H$ if the corresponding stochastic integral
$\int_0^{\cdot}f(B^H_s)dB^H_s$ is of Wick-It\^o type and $f\not\in
C^1({\mathbb R})$. Thus, it is natural to ask whether the similar
It\^o formula for fBm $B^H$ holds or not. This motivates the subject
matters of this paper. Here, we only consider the case
$\frac12<H<1$. This means that fBm $B^H$ is of finite quadratic
variation provided $\frac12<H<1$. Recall that a process $X$ is said
to be of finite quadratic variation if quadratic variation $[X,X]$
is finite. For any continuous finite quadratic variation process $X$
we have (see, for example, Gradinaru {\it et al}~\cite{Grad1})
\begin{equation}\label{sec1-eq100}
F(X_t)=F(0)+\int_0^tf(X_s)d^{-}X_s+\frac12\left[f(X),X\right]_t
\end{equation}
provided $F\in C^2({\mathbb R})$ and $F'=f$, where the integral
$\int_0^tf(X_s)d^{-}X_s$ is the forward (pathwise) integral defined
by
$$
\int_0^tf(X_s)d^{-}X_s=\lim_{\varepsilon\downarrow
0}\frac1{\varepsilon}\int_0^tf(X_s)(X_{s+\varepsilon}-X_s)ds
$$
and $\left[f(X),X\right]_t=\int_0^tf'(X_s)d[X,X]_s$, and more
precisely, for fBm $B^H$ with $\frac12<H<1$ we have
$$
\int_0^tf(B^H_s)d^{-}B^H_s=\int_0^tf(B^H_s)dB^H_s
+H\int_0^tf'(B^H_s)s^{2H-1}ds,
$$
where the integral $\int_0^{\cdot}f(B^H_s)dB^H_s$ is of Wick-It\^o
type. However, the formula~\eqref{sec1-eq100} is only effective on
twice-differentiable functions. It is impossible to list here all
the contributors in previous topics. Some surveys and complete
literatures could be found in Nualart~\cite{Nua4}, Biagini {\it et
al}~\cite{BHOZ}, Hu~\cite{Hu2}, Mishura~\cite{Mishura2},
Russo-Vallois~\cite{Russo-Vallois3} and Gradinaru {\it et
al}~\cite{Grad2,Grad1}.

In the present paper, our aim is devoted to find a substitution tool
of the quadratic variation such that It\^o's formula similar
to~\eqref{sec1-eq1.1} holds for fBm with $\frac12<H<1$ whatever
$f\in C^1({\mathbb R})$. For fBm, in~\cite{Grad1} (see
also~\cite{Grad2} and the references therein) Gradinaru {\it et al}
have introduced some substitution tools and studied some fine
problems. They introduced firstly an It\^o formula with respect to a
symmetric-Stratonovich integral, which is closer to the spirit of
Riemann sums limits, and defined a class of high order integrals
having an interest by themselves. They also treated an It\^o formula
with respect to somehow any symmetric integral, introducing a large
class of symmetric integrals via regularization, and fractional
Brownian motion is not the only process for which their It\^o
formula is valid; there are easy extensions to a more general class
of processes. However, our substitution of quadratic covariation and
technique used here are different from theirs. Moreover, we use the
integral of determinate functions with respect to the weighted local
time of fBm $B^H$
$$
\int_0^t\int_{\mathbb{R}}g(x,s){\mathscr L}^{H}(dx,ds)
$$
and
$$
\int_{\mathbb{R}}f(x){\mathscr L}^{H}(dx,t)
$$
to weaken the hypothesis in some formulas, where ${\mathscr
L}^{H}(x,t)=2H\int_0^t\delta(B^H_s-x)s^{2H-1}ds$ is the weighted
local time of fBm $B^H$. Though our method is only effective on fBm
with $\frac12<H<1$, the merit here has been to concentration fully
on fBm in order to get a stronger statement by fully using fBm's
regularity.

Recently, in~\cite{Grad},~\cite{Nour1} and~\cite{Nour2} Nourdin {\it
et al} studied weighted power variations of fBm and introduced some
fine results. More precisely, they showed that the convergence
$$
n^{2Hp-1}\sum_{k=0}^{n-1}h(B^H_{k/n})
\left(B^H_{(k+1)/n}-B^H_{k/n}\right)^{2p}\longrightarrow
\mu_{2p}\int_0^1h(B^H_s)ds\qquad {\rm a.s.}
$$
holds for all $0<H<1$, as $n$ tends to infinity, where $p\geq 1$ is
an integer number, $h\in C({\mathbb R})$ and $\mu_{2p}$ is the
$2p$-moment of a standard Gaussian random variable $G\sim N(0,1)$.
As a corollary, the following convergence holds:
\begin{equation}\label{sec1-eq1.3}
\sum_{k=0}^{n-1}k^{2Hp-1}h(B^H_{k/n})
\left(B^H_{(k+1)/n}-B^H_{k/n}\right)^{2p}\longrightarrow
\mu_{2p}\int_0^1h(B^H_s)s^{2Hp-1}ds\qquad {\rm a.s.}
\end{equation}
for all $0<H<1$. This inspirits us to consider the {\it weighted
quadratic covariation} $[f(B^H),B^H]^{(W)}$ of $f(B^H)$ and $B^H$,
defined by
\begin{equation*}
\left[f(B^H),B^H\right]^{(W)}_t:=2H\lim_{n\to
\infty}\sum_{k=0}^{n-1}k^{2H-1}
\left\{f(B^H_{t_{k+1}})-f(B^H_{t_k})\right\}
\left(B^H_{t_{k+1}}-B^H_{t_k} \right),
\end{equation*}
where the limit is uniform in probability and $t_k=kt/n$. Clearly,
for $f\in C^1({\mathbb R})$ the convergence~\eqref{sec1-eq1.3}
implies that $[f(B^H),B^H]^{(W)}_t$ exists and equals to
$$
2H\int_0^tf'(B^H_s)s^{2H-1}ds,
$$
and moreover the fractional It\^o formula can be rewritten as
\begin{align*}
F(B^H_t)&=F(0)+\int_0^tf(B^H_s)dB^H_s+H\int_0^tf'(B^H_s)s^{2H-1}ds\\
&\equiv
F(0)+\int_0^tf(B^H_s)dB^H_s+\frac12\left[f(B^H),B^H\right]^{(W)}_t
\end{align*}
for all $t\in [0,T]$ and all $0<H<1$. This is the start point in
this paper and we only consider fBm with $\frac12<H<1$. Our aims are
to study the integrals with respect to the weighted local time of
fBm $B^H$ and the weighted quadratic covariation, and to extend the
above the fractional It\^o formula to $f\not \in C^1({\mathbb R})$.
The case $0<H<\frac12$ will be discussed in a forthcoming paper.

This paper is organized as follows. In Section~\ref{section2} we
present some preliminaries for fBm. To use Young integration to
establish integral with respect to local time, in
Section~\ref{section3}, we first investigate the power variation of
the function $x\mapsto {\mathscr L}^{H}(x,t)$ for every $t\geq 0$.
We show that
\begin{align*}
\lim_{n\to \infty}\sum_{\Delta_n}\left|{\mathscr
L}^{H}(a_{i+1},t)-{\mathscr L}^{H}(a_i,t)\right|^{\frac{2H}{3H-1}}
\end{align*}
exists in $L^1$ for every $0\leq t\leq T$ and $\frac12<H<1$, where
$(\Delta_n)$ is a partition of the interval $[a,b]$ such that
$|\Delta_n|\rightarrow 0$ as $n\rightarrow \infty$. Thus, as a
direct consequence we construct a one parameter integral with
respect to these local times with $\frac12<H<1$. The {\it weighted
quadratic covariation} is considered in Section~\ref{section5}. We
show that if the measurable function $x\mapsto f(x)$ is of bounded
$p$-variation with $1\leq p<\frac{2H}{1-H}$ and $\frac12<H<1$, then
the weighted quadratic covariation $\left[f(B^H),B^H\right]^{(W)}$
exists in $L^1$ and
$$
\left[f(B^H),B^H\right]^{(W)}_t=-\int_{\mathbb R}f(x){\mathscr
L}^{H}(dx,t),\qquad t\in [0,T],
$$
and moreover the It\^o formula
$$
F(B^H_t)=F(0)+\int_0^tf(B^H_s)dB^H_s+\frac12\left[f(B^H),B^H
\right]^{(W)}_t
$$
holds for all absolutely continuous function
$F(x)=F(0)+\int_0^xf(y)dy$. In Section~\ref{sec6} we extend these
results to the time-dependent case.

\section{Preliminaries on fractional Brownian motion}\label{section2}
In this section, we briefly recall some basic definitions and
results of fBm. For more aspects on these materials we refer to
Nualart~\cite{Nua4}, Hu~\cite{Hu2}, Biagini {\it et al}~\cite{BHOZ},
Mishura~\cite{Mishura2} and the references therein. Throughout this
paper we assume that $\frac12<H<1$ is arbitrary but fixed and let
$B^H=\{B_t^H, 0\leq t\leq T\}$ be a one-dimensional fBm with Hurst
index $H$ defined on $(\Omega, \mathcal{F}, P)$. Let $({\mathcal
S})^*$ be the Hida space of stochastic distributions and let
$\diamond$ denote the Wick product on $({\mathcal S})^*$. Then
$t\mapsto B_t^H$ is differentiable in $({\mathcal S})^*$. Denote
$$
W^{(H)}_t=\frac{dB_t^H}{dt}\in ({\mathcal S})^*.
$$
We call $W^{(H)}$ the fractional white noise. For $u: {\mathbb
R}_+\to ({\mathcal S})^*$, in a white noise setting we define its
(generalized) fractional stochastic integral of It\^o type with
respect to $B^H$ by
\begin{equation}\label{eq2.100}
\int_0^tu_sdB^H_s:=\int_0^tu_s\diamond W^{(H)}_sds,
\end{equation}
whenever the last integral exists as an integral in $({\mathcal
S})^*$. We call these fractional It\^o integrals, because these
integrals share many properties of the classical It\^o integral. For
any $F\in C^{2,1}({\mathbb R}\times [0,+\infty))$ satisfying
$$
\int_0^T\left|\frac{\partial^2 F}{\partial^2
s}(B_s^H,s)\right|s^{2H-1}ds<\infty\qquad a.s.,
$$
the It\^o type formula
\begin{align}\tag*{}
F(B_t^H,t)=F(0,0)+&\int_0^t\frac{\partial}{\partial x}
F(B_s^H,s)dB_s^H\\     \label{eq2.1}
&+\int_0^t\frac{\partial}{\partial
s}F(B_s^H,s)ds+H\int_0^t\frac{\partial^2}{\partial
x^2}F(B_s^H,s)s^{2H-1}ds
\end{align}
holds. For the It\^o formula of general Gaussian processes, for
example, see Al\'os {\it et al} ~\cite{Nua1} and
Nualart--Taqqu~\cite{Nua-Taq}.

Recall that fBm $B^H$ has a local time ${\mathcal L}^{H}(x,t)$
continuous in $(x,t)\in {\mathbb R}\times [0,\infty)$ which
satisfies the occupation formula (see Geman-Horowitz~\cite{Geman})
\begin{equation}\label{sec2-1-eq1}
\int_0^t\phi(B_s^{H},s)ds=\int_{\mathbb R}dx\int_0^t\phi(x,s)
{\mathcal L}^{H}(x,ds)
\end{equation}
for every continuous and bounded function $\phi(x,t):{\mathbb
R}\times {\mathbb R}_{+}\rightarrow {\mathbb R}$, and such that
$$
{\mathcal L}^{H}(x,t)=\int_0^t\delta(B_s^H-x)ds=\lim_{\epsilon
\downarrow
0}\frac{1}{2\epsilon}\lambda\big(s\in[0,t],|B_s^H-x|<\epsilon\big),
$$
where $\lambda$ denotes Lebesgue measure and $\delta(x)$ is the
Dirac delta function. Define the so-called weighted local time
${\mathscr L}^{H}(x,t)$ of $B^H$ at $x$ as follows
$$
{\mathscr L}^{H}(x,t)=2H\int_0^ts^{2H-1}{\mathcal L}^{H}(x,ds)\equiv
2H\int_0^t\delta(B_s^H-x)s^{2H-1}ds.
$$
Then the the occupation formula~\eqref{sec2-1-eq1} can be rewritten
as
\begin{equation}
2H\int_0^t\phi(B_s^{H},s)s^{2H-1}ds=\int_{\mathbb
R}dx\int_0^t\phi(x,s) {\mathscr L}^{H}(x,ds),
\end{equation}
and the following Tanaka formula holds:
\begin{equation}\label{eq2.4}
|B_t^H-x|=|x|+\int_0^t{\rm sign}(B^H_s-x)dB^H_s+{\mathscr
L}^{H}(x,t).
\end{equation}

Consider the integral representation of fBm $B^H$ of the form
\begin{equation}\label{sec2-eq100}
B^H_t=\int_0^tK_H(t,u)dB_u,\qquad 0\leq t\leq T,
\end{equation}
where $B$ is a standard Brownian motion and the kernel $K_H(t,u)$
satisfies
\begin{equation}
\frac{\partial K_H}{\partial t}(t,u)=\kappa_H\left(\frac12-H\right)
\left(\frac{u}{t}\right)^{\frac12-H}(t-u)^{H-\frac32}
\end{equation}
with a normalizing constant $\kappa_H>0$ given by
$$
\kappa_H=\left(\frac{2H\Gamma(\frac32-H)}{\Gamma(H+\frac12)
\Gamma(2-2H)}\right)^{1/2}.
$$
Define the operator $\Gamma_{H,T}$ on $L^2([0; T])$ by
$$
\Gamma_{H,T}h(t)=\int_0^t\frac{\partial }{\partial
t}K_H(t,u)h(u)du,\qquad h\in L^2([0; T]),
$$
and let ${\bf S}$ denote the set of all smooth functions on $[0,T]$
with bounded derivatives. Then the function $\Gamma_{H,T}h(t)$ is
continuous and the transpose $\Gamma_{H,t}^*$ of $\Gamma_{H,T}$
restricted to the interval $[0,t]$ ($0\leq t\leq T$) is given by
$$
\Gamma_{H,T}^*g(u)=-\kappa_Hu^{\frac12-H}\frac{d}{du}\int_u^t
m^{H-\frac12}(m-u)^{H-\frac12}g(m)dm,\qquad 0\leq u\leq t
$$
for $g\in {\bf S}$, the set of all smooth functions on $[0,T]$ with
bounded derivatives. In particular, for $\frac12<H<1$ we have
$$
\Gamma_{H,t}^*g(u)=(H-\frac12)\kappa_H u^{\frac12-H}\int_u^t
m^{H-\frac12}(m-u)^{H-\frac32}g(m)dm.
$$
The following lemma gives the computation of
$\mathbf{\Gamma}_{H,t}\mathbf{\Gamma}_{H,t}^{\ast}$ (see~\cite{Hu2},
for the proof).
\begin{lemma}\label{lemma2.1}
Let $f\in {\bf S}$. Then, for $\frac12<H<1$ we have
\begin{equation}\label{sec2-eq8}
\mathbf{\Gamma}_{H,t}\mathbf{\Gamma}_{H,t}^{\ast}f(s)
=H(2H-1)\int_0^t|s-u|^{2H-2}f(u)du.
\end{equation}
\end{lemma}

Finally, recall that fractional Malliavin derivative for smooth
random variable
$$
F=f\left(\int_0^T\eta_1(t)dB_t^H,\ldots,
\int_0^T\eta_n(t)dB_t^H\right),
$$
with $f\in C^{\infty}_b({\mathbb R}^n)$ is defined
\begin{equation}\label{equation1.1}
D_s^HF=\sum_{j=1}^{n}\frac{\partial f}{\partial
x_j}\left(\int_0^T\eta_1(t)dB_t^H,\ldots,\int_0^T
\eta_n(t)dB_t^H\right)\eta_j(s),
\end{equation}
where $\eta_j$ ($j=1,2,\ldots,n$) satisfies
$\int_0^T\left[\Gamma_{H,T}^*\eta_j(t)\right]^2dt <\infty$. Denote $
{\mathbb D}_s^H=\Gamma_{H,T}\Gamma_{H,T}^*D^H_s$. We have
(Proposition 6.24 in~\cite{Hu2})
\begin{align*}
E&\left(\int_0^Tg(s)dB^H_s\right)^2=\int_0^T\left(E\left[\Gamma_{H,T}^*
g(s)\right]\right)^2ds+\int_0^T\int_0^TE\left[{\mathbb
D}_s^Hg(r){\mathbb D}_r^Hg(s)\right]ds\\
&\qquad\equiv \alpha_H\int_0^T\int_0^TE[g(s)g(r)]|s-r|^{2H-2}dsdr
+\int_0^T\int_0^TE\left[{\mathbb D}_s^Hg(r){\mathbb
D}_r^Hg(s)\right]dsdr
\end{align*}
with $\alpha_H=H(2H-1)$, if the right hand side of this identity is
finite.

\section{Power variation of fractional local time}\label{section3}
To use Young integration to establish integral with respect to local
time of fBm, we first investigate $p$-variation of the mapping
$x\mapsto {\mathscr L}^{H}(x,t)$ for every $t\geq 0$.
\begin{definition}\label{definition1.1}
Let $p\geqslant 1$ be a fixed real number. A function
$f:[a,b]\mapsto {\mathbb R}$ is of bounded $p$-variation if
$$
\sup_{\triangle_n}\sum_{i=0}^{n}|f(x_{i+1})-f(x_i)|^p<\infty,
$$
where the supremum is taken over all partition
$\triangle_n=\{a=x_0<x_1<\cdots<x_n=b\}$ of $[a,b]$.
\end{definition}

By the local nondeterminacy of fBm we can prove the following
estimate.
\begin{lemma}\label{lemma3.2}
For all $s,r\in [0,T],\;s\geq r$ and $0<H<1$ we have
\begin{equation}\label{sec3-eq3.2}
s^{2H}r^{2H}-\mu^2\geq \kappa(s-r)^{2H}r^{2H},
\end{equation}
where $\mu=E(B^H_sB^H_r)$ and $\kappa>0$ is a constant.
\end{lemma}

\begin{lemma}\label{lemma3.1}
For $t\geq 0,x\in {\mathbb R}$ set
$$
\widetilde{B}_t^H(x):=\int_0^t1_{(B_s^H>x)}dB_s^H.
$$
Then the estimate
\begin{equation}\label{equation4.1}
E[(\widetilde{B}_t^H(b)-\widetilde{B}_t^H(a))^2]\leqslant
C_{H,t}(b-a)^{1+\alpha}
\end{equation}
holds for all $\frac12<H<1$, $0<\alpha<\frac{2H-1}H$ and $a,b\in
{\mathbb R},a<b$, where $C_{H,t}>0$ is a constant depending only on
$H,t$, so the process $\{\widetilde{B}_t^H(x)\,:\,x\in {\mathbb
R}\}$ has $\alpha$-H\"older continuous paths with $\alpha\in
(0,\frac{3H-1}{2H})$ for every $0\leq t\leq T$.
\end{lemma}
\begin{proof}
For $a,b\in {\mathbb R},a<b$ define the function
$f(x)=1_{(a,b]}(x)$. A straightforward calculation shows that
$$
D_s^H1_{[0,\infty)}(B_t^H)=\delta(B_t^H)1_{[0,t]}(s).
$$
It follows that
\begin{align*}
E\left(\widetilde{B}_t^H(b)-\widetilde{B}_t^H(a)\right)^2
&=E\left(\int_0^tf(B^H_s)dB^H_s\right)^2\\
&=H(2H-1)\int_0^t\int_0^t|s-r|^{2H-2}E1_{(a,b]^2}(B^H_s,B^H_r)dsdr
\end{align*}
for $\frac12<H<1$. Now let us estimate $E1_{(a,b]^2}(B^H_s,B^H_r)$.
We have
\begin{align*}
E1_{(a,b]^2}(B^H_s,B^H_r)&=\int_a^b\int_a^b
\frac{1}{2\pi\rho}\exp\left(-\frac1{2\rho^2}(r^{2H}x^2-2\mu
xy+s^{2H}y^2)\right)dxdy\\
&=\frac1{2\pi}
\int_{\frac{a}{r^{H}}}^{\frac{b}{r^{H}}}e^{-\frac12x^2}dx
\int_{\frac{ar^{H}-\mu x}\rho}^{\frac{br^{H}-\mu
x}\rho}e^{-\frac12y^2}dy\\
&\leq\frac1{\sqrt{2\pi}}\int_{\frac{a}{r^{H}}}^{\frac{b}{r^{H}}}
e^{-\frac12x^2}dx \left(\frac1{\sqrt{2\pi}}\int_{\frac{ar^{H}-\mu
x}\rho}^{\frac{br^{H}-\mu
x}\rho}e^{-\frac12y^2}dy\right)^\alpha\\
&\leq \left(\frac{r^H(b-a)}\rho\right)^\alpha
\int_{\frac{a}{r^{H}}}^{\frac{b}{r^{H}}}e^{-\frac12x^2}dx\\
&\leq\frac{r^{(\alpha-1)H}}{\rho^{\alpha}}(b-a)^{1+\alpha},
\end{align*}
where $\rho=\sqrt{s^{2H}r^{2H}-\mu^2}$. It follows from
Lemma~\ref{lemma3.2} that for $\frac12<H<1$,
\begin{align*}
E(\widetilde{B}_t^H(b)-&\widetilde{B}_t^H(a))^2\leq
4H(2H-1)\kappa(b-a)^{1+\alpha}
\int_0^tds\int_0^s(s-r)^{2H-\alpha H-2}r^{-H}dr\\
&=\frac{4\kappa (2H-1)}{(1-\alpha)}t^{H(1-\alpha)}\left(
\int_0^1(1-x)^{2H-\alpha H-2}x^{-H}dx\right)(b-a)^{1+\alpha}.
\end{align*}
This obtains the estimate~\eqref{equation4.1} for $\frac12<H<1$.
\end{proof}

\begin{theorem}\label{theorem3.1}
Let ${\mathscr L}^{H}(x,t)$ be the weighted local time of fBm $B^H$
with $\frac12<H<1$. Then for every $0\leq t\leq T$, the limit
\begin{align*}
\lim_{n\rightarrow \infty}\sum_{i=0}^{n-1}E\left|{\mathscr
L}^{H}(a_{i+1},t)-{\mathscr L}^{H}(a_i,t)\right|^{\frac{2H}{3H-1}}
\end{align*}
exists, and so
\begin{align*}
\lim_{n\rightarrow \infty}\sum_{i=0}^{n-1}\left|{\mathscr
L}^{H}(a_{i+1},t)-{\mathscr L}^{H}(a_i,t)\right|^{p}=0\qquad {\text
{in $L^1$}}
\end{align*}
for $p>\frac{2H}{3H-1}$, where $\{a_0,a_1,\cdots,a_n\}$ is a
partition of $[a,b]$ such that $\max_i\{|a_{i+1}-a_i|\}\to 0$ as $n$
tends to infinity.
\end{theorem}
\begin{proof}
Keeping the notation in Lemma~\ref{lemma3.1}, we get
$$
{\mathscr L}^{H}(x,t)=2\left(\phi_t(x) -\widetilde{B}_t^H(x)\right).
$$
by fractional Meyer-Tanaka formula~\eqref{eq2.4}, where
\begin{equation}\label{equation2.1}
\phi_t(x)=(B_t^H-x)^{+}-(-x)^{+}.
\end{equation}
It follows that for any $t\geq 0,x\in {\mathbb R}$,
\begin{align*}
\sum_{i=0}^{n-1}&\left|{\mathscr L}^{H}(a_{i+1},t)-{\mathscr
L}^{H}(a_i,t)\right|^{\frac{2H}{3H-1}}\\
&\qquad\leq C_H\left[\sum_{i=0}^{n-1}|\phi_t(a_{i+1})
-\phi_t(a_i)|^{\frac{2H}{3H-1}}+\sum_{i=0}^{n-1}
\left|\widetilde{B}_t^H(a_{i+1})-\widetilde{B}_t^H(a_i)
\right|^{\frac{2H}{3H-1}} \right].
\end{align*}
Noting that the function $\phi_t(x)$ is Lipschitz continuous in $x$
with Lipschitz constant $2$, we get the convergence in $L^1$
$$
\sum_{i=0}^{n-1}\left|\phi_t(a_{i+1})-\phi_t(a_i)
\right|^{\frac{2H}{3H-1}} \leqslant
C_H\sum_{i=0}^{n-1}|a_{i+1}-a_i|^{\frac{2H}{3H-1}}\leqslant
C_H|\Delta n|^{\frac{1-H}{3H-1}}(b-a){\longrightarrow}0,
$$
as $n\to \infty$. We finally have
\begin{align*}
\sum_{i=0}^{n-1}E\left|\widetilde{B}_t^H(a_{i+1})-\widetilde{B}_t^H(a_i)
\right|^{\frac{2H}{3H-1}}\leq C_H(b-a)
\end{align*}
by Lemma~\ref{lemma3.1}. This completes the proof.
\end{proof}

\begin{theorem}\label{theorem3.2}
Let $\frac12<H<1$. Then the weighted local time ${\mathscr
L}^{H}(x,t)$ of fBm $B^H$ is of bounded $p$-variation in $x$ for any
$0\leq t\leq T$, for all $p>\frac{2H}{3H-1}$, almost surely.
\end{theorem}
\begin{proof}
Keeping the notation in the proof of Theorem~\ref{theorem3.1}, we
have
$$
{\mathscr L}^{H}(t,x)=2\big(\phi_t(x)-\widetilde{B}_t^H(x)\big).
$$
Let $[-N,N]$ contain the support of ${\mathscr L}^{H}(t,x)$ in $x$
and let
$$
D:=\{-N=a_0<a_1<\cdots<a_n=N\}
$$
be a partition of $[-N,N]$. Then
$$
\sup_D\sum_i|\phi_t(a_{i+1})-\phi_t(a_i)|^p\leq
\sup_D2^p\sum_i(a_{i+1}-a_i)^p\leqslant 2^p(2N)^p<\infty.
$$
On the other hand, Lemma~\ref{lemma3.1} yields
$$
\left|\widetilde{B}_t^H(x)-\widetilde{B}_t^H(y)\right|\leq
G_{H,T}|x-y|^{\frac{1+m}2},\quad {\rm a.s.}
$$
for all $x,y\in {\mathbb R}$ by Garsia-Rodemich-Rumsey Lemma, where
$G_{H,T}$ is a nonnegative random variable such that
$E(G_{H,T}^q)<\infty$ for all $q\geq 1$, and $0<m<\frac{2H-1}H$. It
follows that for all $p\geq \frac2{1+m}>\frac{2H}{3H-1}$
\begin{align*}
\sup_D\sum_i(\widetilde{B}_t^H(a_{i+1})-\widetilde{B}_t^H(a_i))^p&\leq
G_{H,T}^p\sup_D\sum_i(a_{i+1}-a_i)^{\frac{p(1+m)}2}\\
&=G_{H,T}^p(2N)^{\frac{p(1+m)}2}<\infty,\quad {\rm a.s.}
\end{align*}
Thus, we complete the proof.
\end{proof}
\begin{corollary}\label{cor3.1}
Let $\frac12<H<1$. Then the usual local time ${\mathcal L}^{H}(x,t)$
is of bounded $p$-variation in $x$ for any $0\leq t\leq T$ and all
$p>\frac{2H}{3H-1}$, almost surely.
\end{corollary}

Now we can establish one parameter integral of local times of
fractional Brownian motion. Recall that the Riemann-Stieltjes
integral (see Young~\cite{Young})
$$
\int_a^bf(x)dg(x):=\lim_{|\triangle_n|\rightarrow 0}\sum_{i=0}^{n}
f(\xi_{i+1})(g(x_{i+1})-g(x_i))
$$
exists if $f$ and $g$ have finite $p$-variation and finite
$q$-variation in the interval $[a,b]$ respectively, and $f$ and $g$
have no common discontinuities, where $\xi_{i+1}\in[x_i, x_{i+1}]$,
$p,q\geqslant 1$, $\frac{1}{p}+\frac{1}{q}>1$,
$|\triangle_n|=\max_{0\leqslant i\leqslant n}|x_{i+1}-x_i|.$
\begin{proposition}\label{proposition3.1}
For $\frac12<H<1$, if $f(x)$ is of bounded $p$-variation in $x$ with
$1\leqslant p<\frac{2H}{1-H}$, then the Young integral
$$
\int_a^bf(x){\mathscr L}^{H}(dx,t)
$$
exists in $L^1$ for any $0\leq t\leq T$.
\end{proposition}

Consider the properties and characterization of one parameter
integral. We first have

\begin{align*}
\int_{\mathbb R}f(x){\mathscr L}^{H}(dx,t)&=\lim_{\Delta_n\to
0}\sum_{j=1}^nf(x_{j-1})\left({\mathscr
L}^{H}(x_j,t)-{\mathscr L}^{H}(x_{j-1},t)\right)\\
&=\lim_{\Delta_n\to 0}\left[\sum_{j=1}^nf(x_{j-1}){\mathscr
L}^{H}(x_j,t)-\sum_{j=0}^{n-1}f(x_{j}){\mathscr
L}^{H}(x_{j},t)\right]\\
&=-\lim_{\Delta_n\to
0}\sum_{j=1}^n\left(f(x_{j})-f(x_{j-1})\right){\mathscr
L}^{H}(x_{j},t)
\end{align*}
by adding some points in the partition $\Delta_n$ to make
$$
{\mathscr L}^{H}(x_1,t)=0,\qquad {\mathscr L}^{H}(x_n,t)=0,
$$
which yields the following
\begin{corollary}
Under the conditions of Proposition~\ref{proposition3.1}, the
integral
$$
\int_{\mathbb R}{\mathscr L}^{H}(x,t)df(x)
$$
exists in $L^1$ for any $0\leq t\leq T$ and
$$
\int_{\mathbb R}f(x){\mathscr L}^{H}(dx,t)=-\int_{\mathbb
R}{\mathscr L}^{H}(x,t)df(x).
$$
In particular, for $f\in C^1(\mathbb R)$ we have
\begin{equation}\label{sec4-eq4.1}
\int_{\mathbb R}f(x){\mathscr L}^{H}(dx,t)=-\int_{\mathbb
R}{\mathscr L}^{H}(x,t)f'(x)dx.
\end{equation}
\end{corollary}

Define the mollifier $\rho$ by
\begin{equation}\label{sec4-eq4.2}
\rho(x)=
\begin{cases}
c\exp(\frac{1}{(x-1)^2-1}), & {\text {$x\in(0,2)$}},\\
0, & {\text {$x\not\in(0,2)$}},
\end{cases}
\end{equation}
where $c$ is a normalizing constant such that $\int_{\mathbb
R}\rho(x)dx=1$. Set $\rho_n(x)=n\rho(nx)$. For a locally integrable
function $g(x)$ we define
$$
g_n(x)=\int_{\mathbb R}\rho_n(x-y)g(y)dy,\quad n\geq 1.
$$
Then $g_n(x)$ is smooth and
$$
g_n(x)=\int_0^2\rho(z)g(x-\frac{z}n)dz,\quad n\geq 1.
$$
Then, similar to the proof of Theorem 2.1 in~\cite{Feng} one can
gives the following
\begin{lemma}\label{lem4.1}
Let $g(x)$ be a measurable function with bounded $p$-variation,
where $1\leqslant p<\frac{2H}{1-H}$. Suppose that $g_n(x)$ is
defined as above, then we have
$$
\lim_{n\rightarrow \infty}\int_{\mathbb R}g_n(x){\mathscr
L}^{H}(dx,t)=\int_{\mathbb R}g(x){\mathscr L}^{H}(dx,t)\qquad {\text
{a.s}},
$$
\end{lemma}

Recall that if $F$ is the difference of two convex functions (This
means that $F$ is an absolutely continuous function with derivative
of bounded variation.), then the It\^o-Tanaka formula
\begin{align*}
F(B_t^H)&=F(0)+\int_0^tF_{-}^{'}(B_s^H)dB_s^H+\frac12\int_{\mathbb
R}{\mathscr L}^{H}(x,t)F''(dx)\\
&\equiv F(0)+\int_0^tF_{-}^{'}(B_s^H)dB_s^H-\frac12\int_{\mathbb
R}F_{-}'(x){\mathscr L}^{H}(dx,t)
\end{align*}
holds. This is given by Coutin {\it et al}~\cite{Cout} (see also Hu
{\it et al}~\cite{Hu1}). The following result is a
modification of this formula.
\begin{proposition}\label{th4.1}
Let the measurable function $f:{\mathbb R}\to {\mathbb R}$ be of
bounded $p$-variation with $1\leqslant p<\frac{2H}{1-H}$ and let $F$
be an absolutely continuous function with derivative $F'=f$. If the
integral $\int_0^tf(B^{H}_s)dB^{H}_s$ exists, then
$$
F(B_t^H)=F(0)+\int_0^tf(B_s^H)dB_s^H-\frac12\int_{\mathbb
R}f(x){\mathscr L}^{H}(dx,t).
$$
\end{proposition}
\begin{proof}
For $n\geqslant 1$ we set
$$
F_n(x):=\int_{\mathbb R}\rho_n(x-y)F(y)dy,
$$
where $\rho_n$ is the mollifier defined in~\eqref{sec4-eq4.2}. Then,
it is well-known that for each $x$,
$$
\lim_{n\rightarrow \infty}F_n(x)=F(x),\quad \lim_{n\rightarrow
\infty}F_n'(x)=f(x).
$$
Moreover, $F_n\in C^2(\mathbb R)$ and the fractional It\^{o} formula
yields
\begin{equation}\label{sec4-eq4.3}
F_n(B_t^H)=F_n(0)+\int_0^tF_n'(B_s^H)dB_s^H+
H\int_0^ts^{2H-1}F_n''(B_s^H)ds.
\end{equation}
If $n$ tends to infinity, then it is easy to see that $F_n(B^H_t)$
converges to $F(B^H_t)$ almost surely, and $F_n'(B^H_t)$ converges
to $f(B^H_t)$.

Next we consider the limit of the last term in~\eqref{sec4-eq4.3}.
We have
\begin{align*}
2H\int_0^ts^{2H-1}F_n''(B_s^H)ds &=
\int_{\mathbb R}{\mathscr L}^{H}(x,t)F_n''(x)dx\\
&=-\int_{\mathbb R}F_n'(x){\mathscr L}^{H}(dx,t)\\
&\longrightarrow-\int_{\mathbb R}f(x){\mathscr L}^{H}(dx,t)
\end{align*}
almost surely, as $n$ tends to infinity. Finally, we see that
$$
\int_0^tF_n'(B_s^H)dB_s^H=\int_0^tF_n'(B_s^H)\diamond W^{(H)}_sds\to
\int_0^tf(B_s^H)dB_s^H,\qquad {\text {in $({\mathcal S})^{*}$}}
$$
and also almost surely as $n\to \infty$, by the convergence of the
other terms in~\eqref{sec4-eq4.3}. So, the result follows.
\end{proof}

\section{Weighted quadratic covariation}\label{section5}

In this section, we use one parameter integral
$$
\int_{\mathbb R}f(x){\mathscr L}^{H}(dx,t)
$$
to study the {\it weighted quadratic covariation}
$\left[f(B^H),B^H\right]^{(W)}$. Let $B$ be a standard Brownian
motion. Recall that if $f$ is locally square integrable, then the
quadratic covariation $[f(B),B]$ of $f(B)$ and $B$ exists in $L^1$
and
$$
[f(B),B]_t=-\int_{\mathbb R}f(x){\mathscr L}^{B}(dx,t),
$$
where ${\mathscr L}^{B}(x,t)$ is the local time of $B$ and
$[f(B),B]$ is defined by
$$
[f(B),B]_t=\lim_{n\to
\infty}\sum_{k=0}^{n-1}\{f(B_{(k+1)t/n})-f(B_{kt/n})\}(B_{(k+1)t/n}
-B_{kt/n})
$$
with the limit being uniform in probability. For this see F\"ollmer
{\it et al}~\cite{Follmer} and Eisenbaum~\cite{Eisen1}. However,
this is not true for fBm $B^H$, i.e. in general
$$
\left[f(B^H),B^H\right]_t\neq-\int_{\mathbb R}f(x){\mathscr
L}^{H}(dx,t),\qquad 0\leq t\leq T,
$$
because $[f(B^H),B^H]_t=0$ for $1/2<H<1$. But, as we pointed out
before from Gradinaru--Nourdin~\cite{Grad}, Gradinaru {\it et
al}~\cite{Grad1} and Mishura--Valkeila~\cite{Mishura} one can drive
an inspiration to construct the {\it weighted quadratic covariation}
of $f(B^H)$ and $B^H$.
\begin{definition}
Define the process $\left[f(B^H),B^H\right]^{(W)}_t,\;0\leq t\leq T$
by
\begin{align*}
&\left[f(B^H),B^H\right]^{(W)}_t\\
&\qquad :=2H\lim_{n\to \infty}\sum_{k=0}^{n-1}
k^{2H-1}\{f(B^H_{t(k+1)/n})-f(B^H_{tk/n})\}
(B^H_{t(k+1)/n}-B^H_{tk/n})
\end{align*}
with the limit being uniform in probability. This process is called
the weighted quadratic covariation of $f(B^H)$ and $B^H$.
\end{definition}
Clearly, if $H=\frac12$ the weighted quadratic covariation coincides
with the usual quadratic covariation of Brownian motion $B$, and in
particular we have
$$
\left[B^H,B^H\right]^{(W)}_t=t^{2H}.
$$
A classical result (see Klein and Gin\'e~\cite{Klein} when $p=1$)
shows that, for all $t\in [0,T]$
\begin{align*}
A_n(T,t):=2H\sum_{k=0}^{n-1}
k^{2H-1}&1_{[0,t]}({kT/n})(B^H_{(k+1)T/n}-B^H_{kT/n})^{2p}\\
&=2H\sum_{k=0}^{[n\frac{t}T]-1}
k^{2H-1}\left(B^H_{kt/(n\frac{t}T)}-B^H_{kt/(n\frac{t}T)}
\right)^{2p}\\
&\longrightarrow \mu_{2p}t^{2Hp}, \quad {\rm a.s.}
\end{align*}
as $n$ tends to infinity, and moreover, $A_n(T,t)$ also converges to
$\mu_{2p}t^{2Hp}$ in $L^2$, where $\alpha_{2p}$ is the $2p$-moment
of a standard Gaussian random variable $G\sim N(0,1)$. Thus, the
following theorem is also due to Gradinaru--Nourdin~\cite{Grad}.

\begin{theorem}\label{lem5.1}
Let $g\in C({\mathbb R})$ and let $p \geq 1$ a integer. Then the
convergence
\begin{equation}\label{sec5-eq5.2-1}
\sum_{k=0}^{n-1} k^{2Hp-1}g(B^H_{kt/n})
\left(B^H_{(k+1)t/n}-B^H_{kt/n}\right)^{2p}\longrightarrow
\mu_{2p}\int_0^tg(B^H_s)s^{2Hp-1}ds,\quad {\rm a.s.}
\end{equation}
holds for all $0<H<1$, as $n$ tends to infinity.
\end{theorem}
As a direct consequence of the above theorem, for all $f\in
C^{1}({\mathbb R})$ we see that
\begin{align*}
\lim_{n\to \infty}\sum_{k=0}^{n-1}
k^{2H-1}\{f(B^H_{t(k+1)/n})-f(B^H_{tk/n})\}
(B^H_{t(k+1)/n}-B^H_{tk/n})\longrightarrow\int_0^tf'(B^H_s)s^{2H-1}ds
\end{align*}
almost surely, as $n\to \infty$, and
\begin{equation}\label{sec5-eq5.2-1-}
\left[f(B^H),B^H\right]^{(W)}_t=-\int_{\mathbb R}f(x){\mathscr
L}^{H}(dx,t)
\end{equation}
by occupation formula and the identity~\eqref{sec4-eq4.1}. More
generally, we have
\begin{theorem}\label{th5.2}
Let the measurable function $x\mapsto f(x)$ be of bounded
$p$-variation with $1\leqslant p<\frac{2H}{1-H}$. Then the weighted
quadratic covariation $\left[f(B^H),B^H\right]^{(W)}$ exists, and we
have
\begin{equation}\label{th5.2-eq}
\left[f(B^H),B^H\right]^{(W)}_t=-\int_{\mathbb R}f(x){\mathscr
L}^{H}(dx,t)
\end{equation}
for all $t\in [0,T]$.
\end{theorem}
\begin{proof}
If $f\in C^1({\mathbb R})$, this is~\eqref{sec5-eq5.2-1-}. Let now
$f\not\in C^1({\mathbb R})$. For $m\geqslant 1$ we set
$$
f_m(x):=\int_{\mathbb R}\rho_m(x-y)f(y)dy,
$$
where $\rho_m$ is the mollifier defined as in~\eqref{sec4-eq4.2}.
Then $f_m\in C^1({\mathbb R})$ and we have
$$
\left[f_m(B^H),B^H\right]^{(W)}_t=-\int_{\mathbb R}f_m(x){\mathscr
L}^{H}(dx,t)
$$
for all $t\in [0,T]$. Moreover using Lebesgue's dominated
convergence theorem, one can prove that as $m\to \infty$, for each
$x$, $f_m(x)\to f(x)$. In order to finish this proof we consider the
double sequence
$$
\alpha_{mn}(t):=\sum_{k=0}^{n-1}k^{2H-1}\left\{
f_m(B^H_{(k+1)t/n})-f_m(B^H_{kt/n})
\right\}(B^H_{(k+1)t/n}-B^H_{kt/n}), \quad m,n\geq 1.
$$
Elementary analysis show that the following convergence hold almost
surely:
$$
\lim_{n\to \infty}\alpha_{mn}(t)=\left[f_m(B^{H}),B^{H}
\right]_t^{(W)}=-\int_{\mathbb R}f_m(x){\mathscr L}^{H}(dx,t),\quad
m\geq 1,
$$
$$
\lim_{m\to \infty}\left[f_m(B^{H}),B^{H}
\right]_t^{(W)}=-\int_{\mathbb R}f(x){\mathscr L}^{H}(dx,t),
$$
$$
\lim_{m\to \infty}\alpha_{mn}(t)=\sum_{k=0}^{n-1}k^{2H-1}\left\{
f(B^H_{(k+1)t/n})-f(B^H_{kt/n})
\right\}(B^H_{(k+1)t/n}-B^H_{kt/n}),\quad n\geq 1,
$$
where the first convergence above is uniform in $m$. These imply
that
$$
\lim_{m\to \infty}\lim_{n\to \infty}\alpha_{mn}(t)=\lim_{n\to
\infty}\lim_{m\to \infty}\alpha_{mn}(t)
$$
for all $t\geq 0$, and theorem follows.
\end{proof}
According to Proposition~\ref{th4.1}, we get an extension of It\^o's
formula.
\begin{corollary}\label{cor5.1}
Let the measurable function $f:{\mathbb R}\to {\mathbb R}$ be of
bounded $p$-variation with $1\leqslant p<\frac{2H}{1-H}$ and let $F$
be an absolutely continuous function with the derivative $F'=f$. If
the integral $\int_0^tf(B^{H}_s)dB^{H}_s$ exists, then we have
\begin{equation}\label{se4-eq4.111}
F(B^H_t)=F(0)+\int_0^tf(B^H_s)dB^H_s+\frac12
\left[f(B^H),B^H\right]^{(W)}_t
\end{equation}
for all $t\in [0,T]$.
\end{corollary}
Clearly, this is an analogue of F\"ollmer-Protter-Shiryayev's
formula (see~\cite{Follmer}). It is an improvement in terms of the
hypothesis on $f$ and it is also quite interesting itself. Recall
that a process $X$ is called a finite quadratic variation process if
its quadratic variation $[X,X]$
 exists. FBm with Hurst index $H\geq \frac12$ is of finite
quadratic variation. If $X$ is a finite quadratic variation process
and if $F\in C^2({\mathbb R})$, then the following It\^o's formula
holds:
\begin{equation}\label{sec5-eq4.10000}
F(X_t)=F(0)+\int_0^tf(X_s)d^{-}X_s+\frac12\left[f(X),X\right]_t,
\end{equation}
where $f=F'$ and the integral $\int_0^tf(X_s)d^{-}X_s$ is the
forward (pathwise) integral defined by
$$
\int_0^tf(X_s)d^{-}X_s=\lim_{\varepsilon\downarrow
0}\frac1{\varepsilon}\int_0^tf(X_s)(X_{s+\varepsilon}-X_s)ds
$$
and $\left[f(X),X\right]_t=\int_0^tf'(X_s)d[X,X]_s$. If $X$ is a
continuous semimartingale, then the integral
$$
\int_0^tf(X_s)d^{-}X_s=\int_0^tf(X_s)dX_s
$$
is classical It\^o's integral. For fBm $B^H$ with $\frac12<H<1$ we
have
$$
\int_0^tf(B^H_s)d^{-}B^H_s=\int_0^tf(B^H_s)dB^H_s
+H\int_0^tf'(B^H_s)s^{2H-1}ds.
$$
However, the formula~\eqref{sec5-eq4.10000} is only effective on
twice-differentiable functions. Our results is a modification of
this formula, but our It\^o formula is only effective on fBm with
$\frac12<H<1$. Thus, it is natural to ask whether the It\^o
formula~\eqref{sec5-eq4.10000} for any absolutely continuous
function
$$
F(x)=F(0)+\int_0^xf(y)dy
$$
holds or not. This is an interesting question. We refer to Russo and
Vallois~\cite{Russo-Vallois3} (see also Gradinaru {\it et
al}~\cite{Grad1,Grad2}) for a complete description of stochastic
calculus with respect to finite quadratic variation process.

At the end of this section, we consider the process $Y^H$ of the
form
$$
Y_t^{(H)}:=\int_0^tf(B^H_s)dB^H_s+\left[f(B^H),B^H\right]^{(W)}_t,\qquad
0\leq t\leq T.
$$
Clearly, when $H=\frac12$, the process $Y^{(1/2)}$ is the backward
integral $\int_0^tf(B_s)d^{*}B_s$ with respect to Brownian motion
$B=\{B_t,0\leq t\leq 1\}$, defined by
$$
\int_0^tf(B_s)d^{*}B_s:=\lim_{\Delta_n\to
0}\sum_{0<t_1<t_2<\cdots<t_n<t}f(B_{t_{i+1}})(B_{t_{i+1}}-B_{t_{i}})
$$
in probability, where $\Delta_n=\max_{i}\{|t_{i}-t_{i-1}|\}$.
F\"ollmer-Protter-Shiryayev have proved in~\cite{Follmer} (see also
Eisenbaum~\cite{Eisen1}) that if $f$ is locally square integrable,
then the above limit in probability exists, and
$$
\int_0^tf(B_s)d^{*}B_s=-\int_{1-t}^1f(\widehat{B}_s)d\widehat{B}_s
$$
for $t\in [0,1]$, where $\widehat{B}_t=B_{1-t}$ is the time reversal
on $[0,1]$ of the Brownian motion $B$.

\begin{theorem}\label{th5.3}
Let the measurable function $x\mapsto f(x)$ be of bounded
$p$-variation with $1\leqslant p<\frac{2H}{1-H}$. If the integral
$\int_0^tf(B^H_s)dB^H_s$ exists, then we have
\begin{align*}
-\int_{T-t}^Tf(\widehat{B^H}_s)d\widehat{B^H}_s
=\int_0^tf(B^H_s)dB^H_s+\left[f(B^H),B^H\right]^{(W)}_t
\end{align*}
for all $t\in [0,T]$, and so,
$$
\int_{\mathbb R}f(x){\mathscr L}^{H}(dx,t)=\int_0^tf(B^H_s)dB^H_s
+\int_{T-t}^Tf(\widehat{B^H}_s)d\widehat{B^H}_s,
$$
where $\widehat{B^H}_t=B^H_{T-t}$ is the time reversal on $[0,T]$ of
fBm $B^H$.
\end{theorem}
\begin{proof}
If $f\in C^1({\mathbb R})$, then the It\^o formula for $F$ with
$F'=f$ yields
$$
F(\widehat{B^H}_t)=F(\widehat{B^H}_0)+\int_0^t
f(\widehat{B^H}_s)d\widehat{B^H}_s+H\int_0^tf'(\widehat{B^H}_s)
(T-s)^{2H-1}ds
$$
for all $t\in [0,T]$. Elementary calculus yields
\begin{align*}
\int_{T-t}^Tf(\widehat{B^H}_s)d\widehat{B^H}_s&
=F(0)-F(B^H_t)-H\int_{T-t}^Tf'(\widehat{B^H}_s)(T-s)^{2H-1}ds\\
&=-\int_0^tf(B^H_s)dB^H_s-\left[f(B^H),B^H\right]^{(W)}_t.
\end{align*}
For more details about the It\^o formula of general Gaussian
processes, see Al\'os {\it et al} ~\cite{Nua1} and
Nualart--Taqqu~\cite{Nua-Taq}.

Now, for $f\not\in C^1({\mathbb R})$ we set
$$
f_n(x)=\int_{\mathbb R}\rho_n(x-y)f(y)dy,\quad n\geq 1,
$$
where $\rho_n(x)=n\rho(nx)$ and $\rho(x)$ is the mollifier defined
in~\eqref{sec4-eq4.2}. Then $f_n(x)$ is smooth and
$$
f_n(x)\longrightarrow f(x)
$$
for all $x$, as $n$ tends to infinity, and
\begin{align*}
-\int_{T-t}^Tf_n(\widehat{B^H}_s)d\widehat{B^H}_s
&=\int_0^tf_n(B^H_s)dB^H_s+\left[f_n(B^H),B^H\right]^{(W)}_t\\
&\equiv\int_0^tf_n(B^H_s)dB^H_s-\int_{\mathbb R}f_n(x){\mathscr
L}^H(dx,t)
\end{align*}
for all $n\geq 1$. This completes the proof by taking $n\to \infty$.
\end{proof}

\section{Two parameter integrals of local time}\label{sec6}

Now we turn to the local time-space integral. We will use some idea
from Feng-Zhao~\cite{Feng} to define the two parameter integrals of
local time
$$
\int_0^t\int_{\mathbb R}g(x,s){\mathscr L}^{H}(dx,ds).
$$
Recall that the function $(x,y)\mapsto F(x,y)$, defined on
$[a,b]\times[c,d]$ is of bounded $p$-variation in $x$ uniformly in
$y$, if
$$
\sup_{y\in [c,d]}\sum_{i=1}^m|F(x_i,y)-F(x_{i-1},y)|^p<\infty,
$$
where $\{a=x_0<x_1<\cdots<x_m=b\}$ is an arbitrary partition of
$[a,b]$, and furthermore, it is of bounded $p,q$-variation in
$(x,y)$, if
$$
\sup_{[a,b]\times[c,d] }\sum_{j=1}^n\left(\sum_{i=1}^m|\Delta
F(x_i,y_j)|^p\right)^q<\infty,
$$
where
$$
\Delta F(x_i,y_j)=F(x_i,y_j)-F(x_{i-1},y_j)
-F(x_i,y_{j-1})+F(x_{i-1},y_{j-1}),
$$
and $\{a=x_0<x_1<\cdots<x_m=b;c=y_0<y_1<\cdots<y_n=d\}$ is an
arbitrary partition of $[a,b]\times[c,d]$.

A function $(x,y)\mapsto f(x,y)$ is called to have a jump at
$(x_1,y_1)$ if there exists an $\delta>0$ such that for any
$\varepsilon>0$, there exists $(x_2, y_2)$ satisfying
$\max\{|x_1-x_2|,|y_1-y_2|\}<\delta$ and
$$
\left|f(x_2,y_2)-f(x_1,y_2)-f(x_2,y_1)+f(x_1,y_1)\right|>\varepsilon.
$$
A function $(x,y)\mapsto f(x,y)$ is called to satisfy {\it the
finite large jump condition} if for any $\varepsilon>0$, there
exists at most finite many points $\{x_1,\ldots,x_{n_1}\}$,
$\{y_1,\ldots,y_{m_1}\}$ and a constant $\delta_\varepsilon>0$ such
that the total $(p,q)$-variation of $f$ on
$[x,x+\delta_\varepsilon]\times[y',y'']$ is smaller than
$\varepsilon$ when $[x,x+\delta_\varepsilon]\cap
\{x_1,\ldots,x_{n_1}\}=\emptyset$, and the total $(p,q)$-variation
of $f$ on $[x',x'']\times[y,y+\delta_\varepsilon]$ is smaller than
$\varepsilon$ when $[y,y+\delta_\varepsilon]\cap
\{y_1,\ldots,y_{m_1}\}=\emptyset$.

Let now $G:{\mathbb R}\times{\mathbb R}\to {\mathbb R}$ be a
continuous function of bounded $q_1$-variation in $x$ uniformly in
$y$, and be of bounded $q_2$-variation in $y$ uniformly in $x$;
$F:{\mathbb R}\times{\mathbb R}\to {\mathbb R}$ be of bounded
$p,q$-variation in $(x,y)$ and satisfy the finite large jump
condition, where $p,q,q_1,q_2\geq 1$. Then the Young integral (see
Theorem 3.1 in Feng-Zhao~\cite{Feng})
\begin{align*}
\int_a^b\int_c^dG(x,s)F(dx,dy):=\lim_{\Delta_{m,n}\to
0}\sum_{j=1}^{n}\sum_{i=1}^{m}G(x_{i-1},y_{j-1})\Delta F(x_j,y_j)
\end{align*}
is well-defined, where
$\Delta_{m,n}=\max_{i,j}\{|(x_i,y_j)-(x_{i-1},y_{j-1})|\}$, if there
exist two monotone increasing functions $\rho:{\mathbb R}\to
{\mathbb R}_{+}$ and $\sigma:{\mathbb R}\to {\mathbb R}_{+}$ such
that
$$
\sum_{n,m}\rho\left(\frac1{n^{1/q_1}}\right)
\sigma\left(\frac1{m^{1/q_2}}\right)\frac1{n^{1/p}m^{1/pq}}<\infty.
$$
More works for two-parameter $p,q$-variation path integrals can be
fund in Feng-Zhao~\cite{Feng}. By taking $\rho(u)=u^\alpha$ and
$\sigma(u)=u^{1-\alpha}$ with $\alpha\in (0,1)$, one can prove the
following.

\begin{proposition}\label{pro6.1}
Let the measurable function $F: {\mathbb R}\times [0,t]\to {\mathbb
R}$ be of bounded $p,q$-variation in $(x,t)$, and of bounded
$\gamma$-variation in $x$ uniformly in $t$, and satisfy the finite
large jump condition. If these parameters $p,q,\gamma$ satisfy the
conditions
\begin{equation}\label{sec6-eq6.1}
1\leq \gamma<\frac{2H}{1-H},\qquad 2Hpq<2Hq+3H-1, \quad p,q\geqslant
1,
\end{equation}
then the Young integral of two parameters
\begin{align*}
\int_0^t\int_{\mathbb R}{\mathscr L}^{H}(x,s)F(dx,ds)
\end{align*}
is well-defined, almost surely, and
\begin{align*}
\int_0^t\int_{\mathbb R}F(x,s){\mathscr L}^{H}(dx,ds)
=\int_0^t\int_{\mathbb R}&{\mathscr L}^{H}(x,s)F(dx,ds)\\
&-\int_{\mathbb R}{\mathscr L}^{H}(x,t)F(dx,t).
\end{align*}
\end{proposition}

\begin{corollary}\label{cor6.1}
Under the conditions of Proposition~\ref{pro6.1}, we have
$$
\int_0^t\int_{\mathbb R}F(x,s){\mathscr L}^{H}(dx,ds)=-
\int_{\mathbb R}dx\int_0^t\frac{\partial }{\partial
x}F(x,s){\mathscr L}^{H}(x,ds),
$$
provided $F(x,t)\in C^{1,1}({\mathbb R}\times [0,T])$.
\end{corollary}

\begin{corollary}\label{cor6.2}
Under the conditions of Proposition~\ref{pro6.1} we set
$$
F_m(x,s):=\int_0^2\int_0^2\rho(r)\rho(z)F(x-\frac{r}{m},s-
\frac{z}{m})drdz,\qquad m\geqslant 1,
$$
where $\rho$ is the mollifier defined in~\eqref{sec4-eq4.2}. Then
the convergence
$$
\int_0^t\int_{\mathbb R}F_m(x,s){\mathscr
L}^{H}(dx,ds)\longrightarrow \int_0^t\int_{\mathbb R}F(x,s){\mathscr
L}^{H}(dx,ds)
$$
holds almost surely, as $m\to \infty$.
\end{corollary}

By using Corollary~\ref{cor6.2} and occupation formula, we
immediately get an extension of the It\^{o} formula as follows.
\begin{theorem}\label{th6.1}
Let $F\in C^{1,1}({\mathbb R}\times {\mathbb R}_{+})$. Suppose that
the function $(x,t)\mapsto \frac{\partial }{\partial x}F(x,t)$ is of
bounded $p,q$-variation in $(x,t)$, and of bounded
$\gamma$-variation in $x$ uniformly in $t$, where these parameters
$p,q,\gamma$ satisfy the condition~\eqref{sec6-eq6.1}. If the
integral $\int_0^t\frac{\partial }{\partial x}F(B^H_s,s)dB^H_s$
exists, then the following It\^{o} formula holds:
\begin{align*}
F(B_t^H,t)& =F(0,0)+\int_0^t\frac{\partial }{\partial
t}F(B_s^H,s)ds+\int_0^t \frac{\partial }{\partial x}F(B_s^H,s)dB_s^H\\
&\hspace{3cm}-\frac12\int_0^t\int_{\mathbb R}\frac{\partial
}{\partial x}F(x,s){\mathscr L}^{H}(dx,ds).
\end{align*}
\end{theorem}
\begin{proof}
Let $F_m$ be defined in Corollary~\ref{cor6.2} for $m\geq 1$. Then
$F_m\in C^{2,1}({\mathbb R}\times {\mathbb R}_{+})$ and
$\frac{\partial }{\partial x}F_m$, $m\geq 1$ are of bounded $p,
q$-variation in $(x,t)$. By occupation formula it follows that
\begin{align*}
2H\int_0^t\frac{\partial^2}{\partial
x^2}F_m(B^H_s,s)s^{2H-1}ds&=\int_{\mathbb
R}dx\int_0^t\frac{\partial^2}{\partial x^2}F_m(x,s){\mathscr
L}^{H}(x,ds)\\
&=-\int_0^t\int_{\mathbb R}\frac{\partial }{\partial
x}F_m(x,s){\mathscr L}^{H}(dx,ds)\\
&\to-\int_0^t\int_{\mathbb R}\frac{\partial }{\partial
x}F(x,s){\mathscr L}^{H}(dx,ds),
\end{align*}
almost surly, as $m\to \infty$. On the other hand, the fractional
It\^o formula yields
\begin{align*}
F_m(B_t^H,t)&=F_m(0,0)+\int_0^t\frac{\partial }{\partial
t}F_m(B_s^H,s)ds+\int_0^t \frac{\partial }{\partial x}F_m(B_s^H,s)dB_s^H\\
&\hspace{3cm}+H\int_0^t\frac{\partial^2 }{\partial
x^2}F_m(B^H_s,s)s^{2H-1}ds,
\end{align*}
for all $m\geq 1$. This completes the proof by taking $m\to \infty.$
\end{proof}

As the end of this paper, for a measurable function $f:\;{\mathbb
R}\times [0,T]\to {\mathbb R}$ we consider the {\it weighted
quadratic covariation} $\left[f(B^H,\cdot),B^H\right]^{(W)}$ of
$f(B^H,\cdot)$ and $B^H$, defined by
\begin{align*}
&\left[f(B^H,\cdot),B^H\right]^{(W)}_t\\
&\qquad\quad :=2H\lim_{n\to\infty} \sum_{k=0}^{n-1}
k^{2H-1}\left\{f(B^H_{t_{k+1}},t_{k+1})-f(B^H_{t_{k}},t_{k})\right\}
\left(B^H_{t_{k+1}}-B^H_{t_{k}}\right)
\end{align*}
with the limit being uniform in probability for $t\in [0,T]$, where
$\{t_k:=t\frac{k}{n}, k=0,1,2,\ldots,n\}$ is a special partitions of
$[0,t]$. For Brownian motion $B$ we have known that the quadratic
covariation
$$
[f(B,\cdot),B]_t=\lim_{n\to \infty}\sum_{k=0}^{n-1}
\left\{f(B_{t_{k+1}},t_{k+1})-f(B_{t_{k}},t_{k})\right\}
\left(B_{t_{k+1}}-B_{t_{k}}\right)
$$
exists as a limit uniformly in probability, and
$$
[f(B,\cdot),B]_t=\int_0^tf(B_s,s)d^{*}B_s-\int_0^tf(B_s,s)dB_s,\quad
t\in [0,1],
$$
if $f(\cdot,t)$ is locally square integrable and $f(\cdot,t)$ is
continuous in $t$ as a map from $[0,1]$ to $L^2_{loc}({\mathbb R})$.
For this see Eisenbaum~\cite{Eisen1} and F\"ollmer {\it et
al}~\cite{Follmer}.

\begin{theorem}\label{th6.2}
Let $f\in C({\mathbb R}\times {\mathbb R}_{+})$ be of bounded
$p,q$-variation in $(x,t)$, and of bounded $\gamma$-variation in $x$
uniformly in $t$. If these parameters $p,q,\gamma$ satisfy the
condition~\eqref{sec6-eq6.1}, then we have
\begin{equation}\label{sec6-eq6.2}
\left[f(B^H,\cdot),B^H\right]^{(W)}_t=-\int_0^t\int_{\mathbb
R}f(x,s){\mathscr L}^{H}(dx,ds)
\end{equation}
for all $t\in [0,T]$.
\end{theorem}
\begin{proof}
Let $f\in C^{1,1}({\mathbb R}\times [0,T])$. Thanks to occupation
formula, it is enough to prove
\begin{equation}\label{sec6-eq6.2-1}
\left[f(B^H,\cdot),B^H\right]^{(W)}_t=2H\int_0^t\frac{\partial
}{\partial x}f(B^H_s,s)s^{2H-1}ds
\end{equation}
for all $t\in [0,T]$. Consider the decomposition
\begin{align*}
&\sum_{k=0}^{n-1}
k^{2H-1}\left\{f(B^H_{t_{k+1}},t_{k+1})-f(B^H_{t_{k}},t_{k})
\right\}(B^H_{t_{k+1}}-B^H_{t_{k}})\\
&=\sum_{k=0}^{n-1} k^{2H-1}\frac{\partial }{\partial
x}f(B^H_{t_{k}},t_{k})
(B^H_{t_{k+1}}-B^H_{t_{k}})^2+\sum_{k=0}^{n-1}
k^{2H-1}o(B^H_{t_{k+1}}-B^H_{t_{k}})
(B^H_{t_{k+1}}-B^H_{t_{k}})\\
&\qquad+\frac{t}n\sum_{k=0}^{n-1}k^{2H-1}\frac{\partial }{\partial
t}f(B^H_{t_{k}},t_{k})(B^H_{t_{k+1}}-B^H_{t_{k}})
+\sum_{k=0}^{n-1}k^{2H-1}o(\frac{t}n)(B^H_{t_{k+1}}-B^H_{t_{k}})\\
&\equiv A_1+A_2+A_3+A_4.
\end{align*}
According to Theorem~\ref{lem5.1} we have
$$
A_1=\sum_{k=0}^{n-1} k^{2H-1}\frac{\partial }{\partial
x}f(B^H_{t_{k}},t_{k}) (B^H_{t_{k+1}}-B^H_{t_{k}})^2\longrightarrow
\int_0^t\frac{\partial }{\partial x}f(B^H_s,s)s^{2H-1}ds
$$
almost surely, as $n$ tends to infinity. Moreover, elementary
calculus and the H\"older continuity of fBm can show that
$A_2,A_3,A_4\to 0$ a.s., as n tends to infinity. Thus, we have
obtained~\eqref{sec6-eq6.2-1}.

Let now $f\not\in C^{1,1}({\mathbb R}\times [0,T])$. For $m\geq 1$
we define $f_m$ as follows
$$
f_m(x,s):=\int_0^2\int_0^2\rho(r)\rho(z)f(x-\frac{r}{m},s-
\frac{z}{m})drdz,\qquad m\geqslant 1,
$$
where $\rho$ is the mollifier defined in~\eqref{sec4-eq4.2}. Then
$f_m\in C^{1,1}({\mathbb R}\times {\mathbb R}_{+})$ and we have
$$
\left[f_m(B^H,\cdot),B^H\right]^{(W)}_t=-\int_0^t\int_{\mathbb
R}f_m(x,s){\mathscr L}^{H}(dx,ds)
$$
for all $t\in [0,t]$. By applying the convergence of the double
sequence
$$
\alpha_{mn}(t):=\sum_{k=0}^{n-1}\left\{f_m(B^{H}_{t_{k+1}},t_{k+1})
-f_m(B^{H}_{t_{k}},t_{k})\right\}(B^{H}_{t_{k+1}}-B^{H}_{t_{k}}),
\quad m,n\geq 1,
$$
we can show that $\left[f(B^{H},\cdot),B^{H}\right]^{(W)}_t$ exists
and equals to
$$
-\int_0^t\int_{\mathbb R}f(x,s){\mathscr L}^{H}(dx,ds).
$$
This completes the proof.
\end{proof}
From the proof of Theorem~\ref{th6.2} we see that
$$
\left[f(B^H,\cdot),B^H\right]^{(W)}_t=2H\int_0^t\frac{\partial
}{\partial x}f(B^H_s,s)s^{2H-1}ds
$$
for all $f\in C^{1,1}({\mathbb R}\times {\mathbb R}_{+})$. According
to Theorem~\ref{th6.1} and Theorem~\ref{th6.2}, we get an analogue
of F\"ollmer-Protter-Shiryayev's formula (see~\cite{Follmer}).
\begin{corollary}\label{cor6.3}
Let $F\in C^{1,1}({\mathbb R}\times {\mathbb
R}_{+}),f(x,t)=\frac{\partial }{\partial x}F(x,t)$ and let the
conditions in Theorem~\ref{th6.1} hold. Then we have
\begin{align*}
F(B^H_t,t)=F(0,0)+&\int_0^tf(B^H_s,s)dB^H_s\\
&+\int_0^t\frac{\partial }{\partial
t}F(B^H_s,s)ds+\frac12\left[f(B^H,\cdot),B^H\right]^{(W)}_t
\end{align*}
for all $t\in [0,T]$.
\end{corollary}

\begin{theorem}\label{th6.3}
Let $f\in C({\mathbb R}\times {\mathbb R}_{+})$ be of bounded
$p,q$-variation in $(x,t)$, and of bounded $\gamma$-variation in $x$
uniformly in $t$. Assume that the integral
$\int_0^tf(B^H_s,s)dB^H_s$ exists. If these parameters $p,q,\gamma$
satisfy the condition~\eqref{sec6-eq6.1}, then we have
\begin{align*}
-\int_{T-t}^Tf(\widehat{B^H}_s,T-s)d\widehat{B^H}_s
=\int_0^tf(B^H_s,s)dB^H_s+\left[f(B^H,\cdot),B^H\right]^{(W)}_t
\end{align*}
for all $t\in [0,T]$, where $\widehat{B^H}_t=B^H_{T-t}$ is the time
reversal on $[0,T]$ of the fBm $B^H$.
\end{theorem}

\begin{acknowledgement}
{ The authors would like to thank an anonymous earnest referee whose
remarks and suggestions greatly improved the presentation of our
paper.}
\end{acknowledgement}

\end{document}